
\catcode'32=9
\magnification=1200

\voffset=1cm

\font\tenpc=cmcsc10

\font\eightrm=cmr8
\font\eighti=cmmi8
\font\eightsy=cmsy8
\font\eightbf=cmbx8
\font\eighttt=cmtt8
\font\eightit=cmti8
\font\eightsl=cmsl8
\font\sixrm=cmr6
\font\sixi=cmmi6
\font\sixsy=cmsy6
\font\sixbf=cmbx6

\skewchar\eighti='177 \skewchar\sixi='177
\skewchar\eightsy='60 \skewchar\sixsy='60

\font\tengoth=eufm10
\font\tenbboard=msbm10
\font\eightgoth=eufm7 at 8pt
\font\eightbboard=msbm7 at 8pt
\font\sevengoth=eufm7
\font\sevenbboard=msbm7
\font\sixgoth=eufm5 at 6 pt
\font\fivegoth=eufm5

\font\tengoth=eufm10
\font\tenbboard=msbm10
\font\eightgoth=eufm7 at 8pt
\font\eightbboard=msbm7 at 8pt
\font\sevengoth=eufm7
\font\sevenbboard=msbm7
\font\sixgoth=eufm5 at 6 pt
\font\fivegoth=eufm5

\newfam\gothfam
\newfam\bboardfam

\catcode`\@=11

\def\raggedbottom{\topskip 10pt plus 36pt
\r@ggedbottomtrue}
\def\pc#1#2|{{\bigf@ntpc #1\penalty
\@MM\hskip\z@skip\smallf@ntpc #2}}

\def\tenpoint{%
  \textfont0=\tenrm \scriptfont0=\sevenrm \scriptscriptfont0=\fiverm
  \def\rm{\fam\z@\tenrm}%
  \textfont1=\teni \scriptfont1=\seveni \scriptscriptfont1=\fivei
  \def\oldstyle{\fam\@ne\teni}%
  \textfont2=\tensy \scriptfont2=\sevensy \scriptscriptfont2=\fivesy
  \textfont\gothfam=\tengoth \scriptfont\gothfam=\sevengoth
  \scriptscriptfont\gothfam=\fivegoth
  \def\goth{\fam\gothfam\tengoth}%
  \textfont\bboardfam=\tenbboard \scriptfont\bboardfam=\sevenbboard
  \scriptscriptfont\bboardfam=\sevenbboard
  \def\bboard{\fam\bboardfam}%
  \textfont\itfam=\tenit
  \def\it{\fam\itfam\tenit}%
  \textfont\slfam=\tensl
  \def\sl{\fam\slfam\tensl}%
  \textfont\bffam=\tenbf \scriptfont\bffam=\sevenbf
  \scriptscriptfont\bffam=\fivebf
  \def\bf{\fam\bffam\tenbf}%
  \textfont\ttfam=\tentt
  \def\tt{\fam\ttfam\tentt}%
  \abovedisplayskip=12pt plus 3pt minus 9pt
  \abovedisplayshortskip=0pt plus 3pt
  \belowdisplayskip=12pt plus 3pt minus 9pt
  \belowdisplayshortskip=7pt plus 3pt minus 4pt
  \smallskipamount=3pt plus 1pt minus 1pt
  \medskipamount=6pt plus 2pt minus 2pt
  \bigskipamount=12pt plus 4pt minus 4pt
  \normalbaselineskip=12pt
  \setbox\strutbox=\hbox{\vrule height8.5pt depth3.5pt width0pt}%
  \let\bigf@ntpc=\tenrm \let\smallf@ntpc=\sevenrm
  \let\petcap=\tenpc
  \normalbaselines\rm}
\def\eightpoint{%
  \textfont0=\eightrm \scriptfont0=\sixrm \scriptscriptfont0=\fiverm
  \def\rm{\fam\z@\eightrm}%
  \textfont1=\eighti \scriptfont1=\sixi \scriptscriptfont1=\fivei
  \def\oldstyle{\fam\@ne\eighti}%
  \textfont2=\eightsy \scriptfont2=\sixsy \scriptscriptfont2=\fivesy
  \textfont\gothfam=\eightgoth \scriptfont\gothfam=\sixgoth
  \scriptscriptfont\gothfam=\fivegoth
  \def\goth{\fam\gothfam\eightgoth}%
  \textfont\bboardfam=\eightbboard \scriptfont\bboardfam=\sevenbboard
  \scriptscriptfont\bboardfam=\sevenbboard
  \def\bboard{\fam\bboardfam}%
  \textfont\itfam=\eightit
  \def\it{\fam\itfam\eightit}%
  \textfont\slfam=\eightsl
  \def\sl{\fam\slfam\eightsl}%
  \textfont\bffam=\eightbf \scriptfont\bffam=\sixbf
  \scriptscriptfont\bffam=\fivebf
  \def\bf{\fam\bffam\eightbf}%
  \textfont\ttfam=\eighttt
  \def\tt{\fam\ttfam\eighttt}%
  \abovedisplayskip=9pt plus 2pt minus 6pt
  \abovedisplayshortskip=0pt plus 2pt
  \belowdisplayskip=9pt plus 2pt minus 6pt
  \belowdisplayshortskip=5pt plus 2pt minus 3pt
  \smallskipamount=2pt plus 1pt minus 1pt
  \medskipamount=4pt plus 2pt minus 1pt
  \bigskipamount=9pt plus 3pt minus 3pt
  \normalbaselineskip=9pt
  \setbox\strutbox=\hbox{\vrule height7pt depth2pt width0pt}%
  \let\bigf@ntpc=\eightrm \let\smallf@ntpc=\sixrm
  \normalbaselines\rm}

\tenpoint

\frenchspacing


\newif\ifpagetitre
\newtoks\auteurcourant \auteurcourant={\hfil}
\newtoks\titrecourant \titrecourant={\hfil}

\def\appeln@te{}
\def\vfootnote#1{\def\@parameter{#1}\insert\footins\bgroup\eightpoint
  \interlinepenalty\interfootnotelinepenalty
  \splittopskip\ht\strutbox 
  \splitmaxdepth\dp\strutbox \floatingpenalty\@MM
  \leftskip\z@skip \rightskip\z@skip
  \ifx\appeln@te\@parameter\indent \else{\noindent #1\ }\fi
  \footstrut\futurelet\next\fo@t}

\pretolerance=500 \tolerance=1000 \brokenpenalty=5000
\newdimen\hmargehaute \hmargehaute=0cm
\newdimen\lpage \lpage=13.3cm
\newdimen\hpage \hpage=20cm
\newdimen\lmargeext \lmargeext=1cm
\hsize=11.25cm
\vsize=18cm
\parskip 0pt
\parindent=12pt

\def\margehaute{\vbox to \hmargehaute{\vss}}%
\def\margebasse{\vss}

\output{\shipout\vbox to \hpage{\margehaute\nointerlineskip
  \corpsdepage\margebasse}
  \advancepageno \global\pagetitrefalse
  \ifnum\outputpenalty>-20000 \else\dosupereject\fi}

\def\corpsdepage{\hbox to \lpage{\hss\pagetexte\hskip\lmargeext}}
\def\pagetexte{\vbox{\makeheadline\pagebody\makefootline}}
\headline={\ifpagetitre\titleheadline \else
  \ifodd\pageno\rightheadline \else\leftheadline\fi\fi}
\def\leftheadline{\eightpoint\hfil\the\auteurcourant\hfil}
\def\rightheadline{\eightpoint\hfil\the\titrecourant\hfil}
\def\titleheadline{\hfill}
\pagetitretrue

\def\footnoterule{\kern-6\p@
  \hrule width 2truein \kern 5.6\p@} 

\def\pd#1#2 {\pc#1#2| }

\def\pointir{\discretionary{.}{}{.\kern.35em---\kern.7em}\nobreak
\hskip 0em plus .3em minus .4em }

\def\abstract#1{\vbox{\eightpoint \pc ABSTRACT|\pointir #1}}

\def\titre#1|{\message{#1}
              \par\vskip 30pt plus 24pt minus 3pt\penalty -1000
              \vskip 0pt plus -24pt minus 3pt\penalty -1000
              \centerline{\bf #1}
              \vskip 5pt
              \penalty 10000 }

\def\section#1|{\par\vskip .3cm
                {\bf #1}\pointir}

\def\ssection#1|{\par\vskip .2cm
                {\it #1}\pointir}

\long\def\th#1|#2\finth{\par\medskip
              {\petcap #1\pointir}{\it #2}\par\smallskip}

\long\def\tha#1|#2\fintha{\par\medskip
                    {\petcap #1.}\par\nobreak{\it #2}\par\smallskip}

\def\rem#1|{\par\medskip
            {{\it #1}.\quad}}

\def\rema#1|{\par\medskip
             {{\it #1.}\par\nobreak }}

\def\article#1|#2|#3|#4|#5|#6|#7|
    {{\leftskip=7mm\noindent
     \hangindent=2mm\hangafter=1
     \llap{[#1]\hskip.35em}{#2}.\quad
     #3, {\sl #4}, vol.\nobreak\ {\bf #5}, {\oldstyle #6},
     p.\nobreak\ #7.\par}}
\def\livre#1|#2|#3|#4|
    {{\leftskip=7mm\noindent
    \hangindent=2mm\hangafter=1
    \llap{[#1]\hskip.35em}{#2}.\quad
    {\sl #3}.\quad #4.\par}}
\def\divers#1|#2|#3|
    {{\leftskip=7mm\noindent
    \hangindent=2mm\hangafter=1
     \llap{[#1]\hskip.35em}{#2}.\quad
     #3.\par}}
\mathchardef\conj="0365

\def\qed{\quad\raise -2pt\hbox{\vrule\vbox to 10pt{\hrule width 4pt
\vfill\hrule}\vrule}}

\def\decale#1|{\par\noindent\hskip 28pt\llap{#1}\kern 5pt}

\catcode`\@=12


\catcode`\@=11
\def\matrice#1{\null \,\vcenter {\normalbaselines \m@th
\ialign {\hfil $##$\hfil &&\  \hfil $##$\hfil\crcr
\mathstrut \crcr \noalign {\kern -\baselineskip } #1\crcr
\mathstrut \crcr \noalign {\kern -\baselineskip }}}\,}

\def\petitematrice#1{\left(\null\vcenter {\normalbaselines \m@th
\ialign {\hfil $##$\hfil 
&&\thinspace  \hfil $##$\hfil\crcr
\mathstrut \crcr \noalign {\kern -\baselineskip } #1\crcr
\mathstrut \crcr \noalign {\kern -\baselineskip }}}\right)}

\catcode`\@=12

\def\qed{\quad\raise -2pt\hbox{\vrule\vbox to 10pt{\hrule width 4pt
   \vfill\hrule}\vrule}}


%


\def\il{\bigl]\kern-.25em\bigl]}
\def\ir{\bigr]\kern-.25em\bigr]}

\def\iil{\bigl>\kern-.25em\bigl>}
\def\iir{\bigr>\kern-.25em\bigr>}



\def\Grille{\setbox13=\vbox to 5mm{\hrule width 110mm\vfill}
\setbox13=\vbox{\offinterlineskip
   \copy13\copy13\copy13\copy13\copy13\copy13\copy13\copy13
   \copy13\copy13\copy13\copy13\box13\hrule width 110mm}
\setbox14=\hbox to 5mm{\vrule height 65mm\hfill}
\setbox14=\hbox{\copy14\copy14\copy14\copy14\copy14\copy14
   \copy14\copy14\copy14\copy14\copy14\copy14\copy14\copy14
   \copy14\copy14\copy14\copy14\copy14\copy14\copy14\copy14\box14}
\ht14=0pt\dp14=0pt\wd14=0pt
\setbox13=\vbox to 0pt{\vss\box13\offinterlineskip\box14}
\wd13=0pt\box13}


\def\fleche(#1,#2)\dir(#3,#4)\long#5{%
\noalign{\nointerlineskip\leftput(#1,#2){\vector(#3,#4){#5}}\nointerlineskip}}


\def\hfl#1#2#3{\smash{\mathop{\hbox to#3{\rightarrowfill}}\limits
^{\scriptstyle#1}_{\scriptstyle#2}}}

\def\gfl#1#2#3{\smash{\mathop{\hbox to#3{\leftarrowfill}}\limits
^{\scriptstyle#1}_{\scriptstyle#2}}}


 \message{`lline' & `vector' macros from LaTeX}
 \catcode`@=11
\def\{{\relax\ifmmode\lbrace\else$\lbrace$\fi}
\def\}{\relax\ifmmode\rbrace\else$\rbrace$\fi}
\def\newcount{\alloc@0\count\countdef\insc@unt}
\def\newdimen{\alloc@1\dimen\dimendef\insc@unt}
\def\newwrite{\alloc@7\write\chardef\sixt@@n}

\newwrite\@unused
\newcount\@tempcnta
\newcount\@tempcntb
\newdimen\@tempdima
\newdimen\@tempdimb
\newbox\@tempboxa

\def\@spaces{\space\space\space\space}
\def\@whilenoop#1{}
\def\@whiledim#1\do #2{\ifdim #1\relax#2\@iwhiledim{#1\relax#2}\fi}
\def\@iwhiledim#1{\ifdim #1\let\@nextwhile=\@iwhiledim
        \else\let\@nextwhile=\@whilenoop\fi\@nextwhile{#1}}
\def\@badlinearg{\@latexerr{Bad \string\line\space or \string\vector
   \space argument}}
\def\@latexerr#1#2{\begingroup
\edef\@tempc{#2}\expandafter\errhelp\expandafter{\@tempc}%
\def\@eha{Your command was ignored.
^^JType \space I <command> <return> \space to replace it
  with another command,^^Jor \space <return> \space to continue without it.}
\def\@ehb{You've lost some text. \space \@ehc}
\def\@ehc{Try typing \space <return>
  \space to proceed.^^JIf that doesn't work, type \space X <return> \space to
  quit.}
\def\@ehd{You're in trouble here.  \space\@ehc}

\typeout{LaTeX error. \space See LaTeX manual for explanation.^^J
 \space\@spaces\@spaces\@spaces Type \space H <return> \space for
 immediate help.}\errmessage{#1}\endgroup}
\def\typeout#1{{\let\protect\string\immediate\write\@unused{#1}}}

\font\tenln    = line10
\font\tenlnw   = linew10

\newdimen\@wholewidth
\newdimen\@halfwidth
\newdimen\unitlength 

\unitlength =1pt


\def\thinlines{\let\@linefnt\tenln \let\@circlefnt\tencirc
  \@wholewidth\fontdimen8\tenln \@halfwidth .5\@wholewidth}
\def\thicklines{\let\@linefnt\tenlnw \let\@circlefnt\tencircw
  \@wholewidth\fontdimen8\tenlnw \@halfwidth .5\@wholewidth}

\def\linethickness#1{\@wholewidth #1\relax \@halfwidth .5\@wholewidth}

\newif\if@negarg

\def\lline(#1,#2)#3{\@xarg #1\relax \@yarg #2\relax
\@linelen=#3\unitlength
\ifnum\@xarg =0 \@vline
  \else \ifnum\@yarg =0 \@hline \else \@sline\fi
\fi}

\def\@sline{\ifnum\@xarg< 0 \@negargtrue \@xarg -\@xarg \@yyarg -\@yarg
  \else \@negargfalse \@yyarg \@yarg \fi
\ifnum \@yyarg >0 \@tempcnta\@yyarg \else \@tempcnta -\@yyarg \fi
\ifnum\@tempcnta>6 \@badlinearg\@tempcnta0 \fi
\setbox\@linechar\hbox{\@linefnt\@getlinechar(\@xarg,\@yyarg)}%
\ifnum \@yarg >0 \let\@upordown\raise \@clnht\z@
   \else\let\@upordown\lower \@clnht \ht\@linechar\fi
\@clnwd=\wd\@linechar
\if@negarg \hskip -\wd\@linechar \def\@tempa{\hskip -2\wd\@linechar}\else
     \let\@tempa\relax \fi
\@whiledim \@clnwd <\@linelen \do
  {\@upordown\@clnht\copy\@linechar
   \@tempa
   \advance\@clnht \ht\@linechar
   \advance\@clnwd \wd\@linechar}%
\advance\@clnht -\ht\@linechar
\advance\@clnwd -\wd\@linechar
\@tempdima\@linelen\advance\@tempdima -\@clnwd
\@tempdimb\@tempdima\advance\@tempdimb -\wd\@linechar
\if@negarg \hskip -\@tempdimb \else \hskip \@tempdimb \fi
\multiply\@tempdima \@m
\@tempcnta \@tempdima \@tempdima \wd\@linechar \divide\@tempcnta \@tempdima
\@tempdima \ht\@linechar \multiply\@tempdima \@tempcnta
\divide\@tempdima \@m
\advance\@clnht \@tempdima
\ifdim \@linelen <\wd\@linechar
   \hskip \wd\@linechar
  \else\@upordown\@clnht\copy\@linechar\fi}

\def\@hline{\ifnum \@xarg <0 \hskip -\@linelen \fi
\vrule height \@halfwidth depth \@halfwidth width \@linelen
\ifnum \@xarg <0 \hskip -\@linelen \fi}

\def\@getlinechar(#1,#2){\@tempcnta#1\relax\multiply\@tempcnta 8
\advance\@tempcnta -9 \ifnum #2>0 \advance\@tempcnta #2\relax\else
\advance\@tempcnta -#2\relax\advance\@tempcnta 64 \fi
\char\@tempcnta}

\def\vector(#1,#2)#3{\@xarg #1\relax \@yarg #2\relax
\@linelen=#3\unitlength
\ifnum\@xarg =0 \@vvector
  \else \ifnum\@yarg =0 \@hvector \else \@svector\fi
\fi}

\def\@hvector{\@hline\hbox to 0pt{\@linefnt
\ifnum \@xarg <0 \@getlarrow(1,0)\hss\else
    \hss\@getrarrow(1,0)\fi}}

\def\@vvector{\ifnum \@yarg <0 \@downvector \else \@upvector \fi}

\def\@svector{\@sline
\@tempcnta\@yarg \ifnum\@tempcnta <0 \@tempcnta=-\@tempcnta\fi
\ifnum\@tempcnta <5
  \hskip -\wd\@linechar
  \@upordown\@clnht \hbox{\@linefnt  \if@negarg
  \@getlarrow(\@xarg,\@yyarg) \else \@getrarrow(\@xarg,\@yyarg) \fi}%
\else\@badlinearg\fi}

\def\@getlarrow(#1,#2){\ifnum #2 =\z@ \@tempcnta='33\else
\@tempcnta=#1\relax\multiply\@tempcnta \sixt@@n \advance\@tempcnta
-9 \@tempcntb=#2\relax\multiply\@tempcntb \tw@
\ifnum \@tempcntb >0 \advance\@tempcnta \@tempcntb\relax
\else\advance\@tempcnta -\@tempcntb\advance\@tempcnta 64
\fi\fi\char\@tempcnta}

\def\@getrarrow(#1,#2){\@tempcntb=#2\relax
\ifnum\@tempcntb < 0 \@tempcntb=-\@tempcntb\relax\fi
\ifcase \@tempcntb\relax \@tempcnta='55 \or
\ifnum #1<3 \@tempcnta=#1\relax\multiply\@tempcnta
24 \advance\@tempcnta -6 \else \ifnum #1=3 \@tempcnta=49
\else\@tempcnta=58 \fi\fi\or
\ifnum #1<3 \@tempcnta=#1\relax\multiply\@tempcnta
24 \advance\@tempcnta -3 \else \@tempcnta=51\fi\or
\@tempcnta=#1\relax\multiply\@tempcnta
\sixt@@n \advance\@tempcnta -\tw@ \else
\@tempcnta=#1\relax\multiply\@tempcnta
\sixt@@n \advance\@tempcnta 7 \fi\ifnum #2<0 \advance\@tempcnta 64 \fi
\char\@tempcnta}

\def\@vline{\ifnum \@yarg <0 \@downline \else \@upline\fi}

\def\@upline{\hbox to \z@{\hskip -\@halfwidth \vrule
  width \@wholewidth height \@linelen depth \z@\hss}}

\def\@downline{\hbox to \z@{\hskip -\@halfwidth \vrule
  width \@wholewidth height \z@ depth \@linelen \hss}}

\def\@upvector{\@upline\setbox\@tempboxa\hbox{\@linefnt\char'66}\raise
     \@linelen \hbox to\z@{\lower \ht\@tempboxa\box\@tempboxa\hss}}

\def\@downvector{\@downline\lower \@linelen
      \hbox to \z@{\@linefnt\char'77\hss}}

\thinlines

\newcount\@xarg
\newcount\@yarg
\newcount\@yyarg
\newcount\@multicnt
\newdimen\@xdim
\newdimen\@ydim
\newbox\@linechar
\newdimen\@linelen
\newdimen\@clnwd
\newdimen\@clnht
\newdimen\@dashdim
\newbox\@dashbox
\newcount\@dashcnt
 \catcode`@=12


\newbox\tbox
\newbox\tboxa

\def\leftzer#1{\setbox\tbox=\hbox to 0pt{#1\hss}%
     \ht\tbox=0pt \dp\tbox=0pt \box\tbox}

\def\rightzer#1{\setbox\tbox=\hbox to 0pt{\hss #1}%
     \ht\tbox=0pt \dp\tbox=0pt \box\tbox}

\def\centerzer#1{\setbox\tbox=\hbox to 0pt{\hss #1\hss}%
     \ht\tbox=0pt \dp\tbox=0pt \box\tbox}

%
\def\image(#1,#2)#3{\vbox to #1{\offinterlineskip
    \vss #3 \vskip #2}}


\def\leftput(#1,#2)#3{\setbox\tboxa=\hbox{%
    \kern #1\unitlength
    \raise #2\unitlength\hbox{\leftzer{#3}}}%
    \ht\tboxa=0pt \wd\tboxa=0pt \dp\tboxa=0pt\box\tboxa}

\def\rightput(#1,#2)#3{\setbox\tboxa=\hbox{%
    \kern #1\unitlength
    \raise #2\unitlength\hbox{\rightzer{#3}}}%
    \ht\tboxa=0pt \wd\tboxa=0pt \dp\tboxa=0pt\box\tboxa}

\def\centerput(#1,#2)#3{\setbox\tboxa=\hbox{%
    \kern #1\unitlength
    \raise #2\unitlength\hbox{\centerzer{#3}}}%
    \ht\tboxa=0pt \wd\tboxa=0pt \dp\tboxa=0pt\box\tboxa}

\unitlength=1mm

\def\put(#1,#2)#3{\noalign{\nointerlineskip
                               \centerput(#1,#2){$#3$}
                                \nointerlineskip}}
\def\segment(#1,#2)\dir(#3,#4)\long#5{%
\leftput(#1,#2){\lline(#3,#4){#5}}}
\auteurcourant={DOMINIQUE FOATA {\sevenrm AND} GUO-NIU HAN}
\titrecourant={FINITE DIFFERENCE CALCULUS}
\def\grn{\mathop{\rm grn}\nolimits}
\def\eoc{\mathop{\rm eoc}\nolimits}
\def\pom{\mathop{\rm pom}\nolimits}
\vglue.5cm

\rightline{April 5, 2013}
\bigskip\bigskip
\centerline{\bf Finite Difference Calculus
for Alternating Permutations}
\bigskip
\centerline{\sl Dominique Foata and Guo-Niu Han}
\bigskip
{\narrower\narrower
\eightpoint
\noindent
{\bf Abstract}.\quad
The finite difference equation system introduced by Christiane Poupard in the study of 
tangent trees is 
reinterpreted in the alternating permutation environment.
It makes it possible to make a joint study of both
tangent and secant trees and 
calculate the generating polynomial for alternating permutations 
by a new statistic, referred to as being the greater neighbor of the maximum.

}
\bigskip\medskip
\centerline{\bf 1. Introduction}

\medskip
Let $f=(f_n(k))$ $(n\ge 1,\,1\le k\le 2n-1)$ be a family of rational numbers, displayed in a triangular array of the form
$$
f=\matrice{&&&&f_1(1)\cr
&&&f_2(1)&f_2(2)&f_2(3)\cr
&&f_3(1)&f_3(2)&f_3(3)&f_3(4)&f_3(5)\cr
&f_{4}(1)&f_4(2)&f_4(3)&f_4(4)&f_4(5)&f_{4}(6)&f_{4}(7)\cr
\cdots&\cdots&\cdots&\cdots&\cdots&\cdots&\cdots&\cdots&\cdots\cr
}\leqno(1.1)$$
and consider the
finite difference equation system
$$
\Delta^2 f_n(k)+4\,f_{n-1}(k)
=0\quad(n\ge 2,\,1\le k\le 2n-3),\leqno(1.2)
$$
where $\Delta$ stands for the classical finite difference
operator (see, e.g., [Jo39]) 
$$\leqalignno{
\Delta f_n(k)&:= f_n(k+1)-f_n(k),&(1.3)\cr
\noalign{\hbox{so that}} 
\Delta^2
f_n(k)&=f_n(k+2)-2f_n(k+1)+f_n(k).&(1.4)\cr}
$$

If at each step $n\ge 2$ the two entries $f_{n}(1)$ and
$f_{n}(2)$ are given explicit values, the whole system (1.2) has a {\it unique} solution, as the equation 
$\Delta^2 f_n(1)+4\,f_{n-1}(1)=0$ yields the value of $f_{n}(3)$, then $\Delta^2 f_n(2)+4\,f_{n-1}(2)=0$ the value of $f_{n}(4)$, etc.

The same conclusion holds if the two bordered diagonals
$$\displaylines{\qquad
(f_1(1),f_2(1),f_3(1),f_{4}(1),\ldots,f_{n}(1),\ldots\,),\hfill\cr
\hfill{}
(f_1(1),f_2(3),f_3(5),f_{4}(7),\ldots,f_{n}(2n-1),\ldots\,)\qquad\cr}
$$
are taken as initial values. To see this we first note that the equation
$f_2(1)-2f_2(2)+f_2(3)+4f_1(1)=0$ determines $f_2(2)$ uniquely. Assuming that the triangle $(f_{n'}(m))$ $(1\le m\le 2n'-1,\,n'\le n)$ has been determined, the system $\Delta^2 f_{n+1}(m)+4\,f_{n}(m)=0$ $(1\le m\le 2n-1)$ consists of $(2n-1)$ linear equations with $(2n-1)$ unknowns,  namely, $f_{n+1}(2)$, $f_{n+1}(3)$, \dots~, $f_{n+1}(2n)$, the underlying matrix being trigonal of the form
$$F_{n+1}:=\petitematrice{-2&1&0&0&\cdots&0&0&0\cr
1&-2&1&0&\cdots&0&0&0\cr
0&1&-2&1&\cdots&0&0&0\cr
\vdots&\vdots&\vdots&\vdots&\ddots&\vdots&\vdots&\vdots\cr
0&0&0&0&\cdots&1&-2&1\cr
0&0&0&0&\cdots&0&1&-2\cr
}.
$$
As $\det F_{n+1}=-2n$ $(n\ge 1)$, the system
has a unique solution.

The purpose of this paper is to solve (1.2) in four cases, when the sets of initial values
called {\tt [tan1]}, {\tt [tan2]}, {\tt [sec1]},
{\tt [sec2]} are the following:

\smallskip
\noindent
{\tt [tan1]} $f_{1}(1)=1$; $f_{n}(1)=0$ and 
$f_{n}(2)=2\sum\limits_{k}f_{n-1}(k)$ for $n\ge 2$;

\noindent
{\tt [tan2]} $f_{1}(1)=1$; $f_{n}(1)=f_{n}(2n-1)=0$ 
for $n\ge 2$;

\smallskip
\noindent
{\tt [sec1]} $f_{1}(1)=1$; $f_{n}(1)=\sum\limits_{k}f_{n-1}(k)$ and
$f_{n}(2)=3\sum\limits_{k}f_{n-1}(k)$ for $n\ge 2$;

\noindent
{\tt [sec2]} $f_{1}(1)=1$;
$f_{n}(1)=f_{n}(2n-1)=\sum\limits_{k}f_{n-1}(k)$
for $n\ge 2$;

It will be proved (see Theorem~1.5) that both initial values {\tt [tan1]} and {\tt [tan2]}
(resp. {\tt [sec1]} and {\tt [sec2]}) in fact lead to the same solution of the system and, furthermore, that the solutions found for the $f_n(k)$'s are 
non-negative {\it integral values}.
To avoid any confusion the solutions of (1.2) will be denoted by $(g_{n}(k))$ (resp. $(h_{n}(k))$) when using {\tt [tan1]}
(resp. {\tt [sec1]}). The first numerical values of those solutions are displayed in Fig.~1.1.

{\eightpoint

$$\displaylines{
\matrice{&&&&&g_{1}(1)\!=\!1&&&&&&\quad T_1\!=\!1\cr
&&&&g_{2}(1)\!=\!0&g_{2}(2)\!=\!2&g_{2}(3)\!=\!0&&&&
&\quad T_3\!=\!2\cr 
&&&g_{3}(1)\!=\!0&g_{3}(2)\!=\!4&g_{3}(3)\!=\!8&
g_{3}(4)\!=\!4&g_{3}(5)\!=\!0&&&&\quad T_5\!=\!16\cr 
&&g_{4}(1)\!=\!0&g_{4}(2)\!=\!32&g_{4}(3)\!=\!64
&g_{4}(4)\!=\!80&
g_{4}(5)\!=\!64&g_{4}(6)\!=\!32&g_{4}(7)\!=\!0&&&\quad 
T_7\!=\!272\cr
}\cr 
\noalign{\medskip}
\matrice{&&&h_{1}(1)\!=\!1&&&&\hskip-10pt E_2\!=\!1\cr
&&h_{2}(1)\!=\!1&h_{2}(2)\!=\!3&h_{2}(3)\!=\!1&&&\hskip-10pt
E_4\!=\!5\cr 
&h_{3}(1)\!=\!5&h_{3}(2)\!=\!15&h_{3}(3)\!=\!21&
h_{3}(4)\!=\!15&
h_{3}(5)\!=\!5&&\hskip-10pt E_6\!=\!61\cr
h_{4}(1)\!=\!61&h_{4}(2)\!=\!183&h_{4}(3)\!=\!285&
h_{4}(4)\!=\!327&h_{4}(5)\!=\!285&h_{4}(6)\!=\!183&
h_{4}(7)\!=\!61&\cr 
&&&&&&&\hskip-20pt [E_8\!\!=\!\!1385}\cr
 }
$$
\vskip-8pt
}
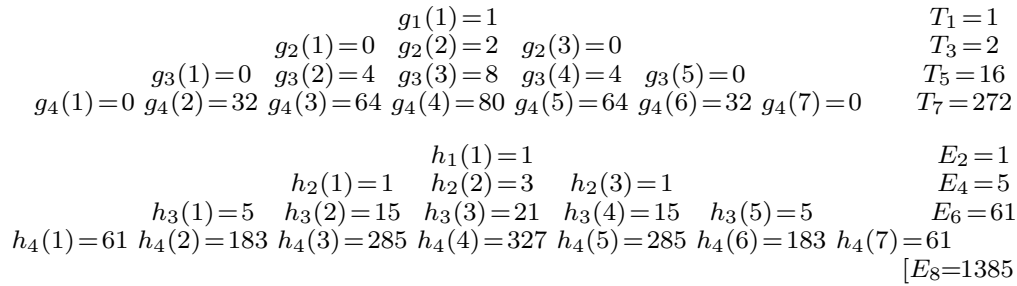
\centerline{Fig. 1.1. The two triangles of the $g_n(k)$'s and $h_{n}(k)$'s.}

\medskip
To the right of each triangle have been calculated the row sums, which are equal, as stated in the next theorem, to the
{\it tangent numbers} (resp. the {\it secant numbers}). Those classical numbers, denoted by $T_{2n+1}$ and $E_{2n}$, appear in the Taylor expansions of $\tan u$ and $\sec u$:
$$
\leqalignno{ \noalign{\vskip-5pt}
\quad\tan u&=\sum_{n\ge 0} {u^{2n+1}\over
(2n+1)!}T_{2n+1}&(1.5)\cr
\noalign{\vskip-5pt}
&={u\over 1!}1+{u^3\over 3!}2+{u^5\over 5!}16+{u^7\over
7!}272+ {u^9\over 9!}7936+\cdots\cr
\qquad\sec u={1\over \cos u}&=\sum_{n\ge 0}
{u^{2n}\over (2n)!}E_{2n}&(1.6)\cr
\noalign{\vskip-5pt}
&=1+{u^2\over 2!}1+{u^4\over 4!}5+{u^6\over 6!}61+{u^8\over
8!}1385+ {u^{10}\over 10!}50521+\cdots\cr 
}
$$
(see, e.g.,
[Ni23, p.~177-178], [Co74, p.~258-259]). 

\proclaim Theorem 1.1. Let $(g_{n}(k))$ (resp.
$(h_{n}(k))$ be the unique solution of the finite difference equation system $(1.2)$ when using the initial values
{\tt [tan1]} (resp. {\tt [sec1]}). Then, the row sums of the solutions are equal to
$$
\leqalignno{\sum_{k}g_{n}(k)&=T_{2n-1}\quad (n\ge 1);&(1.7)\cr
\sum_{k}h_{n}(k)&=E_{2n}\quad (n\ge 1).&(1.8)\cr
}
$$

As further mentioned, Theorem 1.1 will appear as a consequence of Theorem 1.4.
It will also be shown that the generating functions for the coefficients $g_{n}(k)$ and $h_{n}(k)$  can be evaluated in the following forms.

\proclaim Theorem  1.2. Let
$$
Z(x,y):=1+\sum_{n\ge 1}\sum_{1\le k\le 2n+1}
f_{n+1}(k){x^{2n+1-k}\over (2n+1-k)!}
{y^{k-1}\over (k-1)!}
$$
and $Z^{\tan}(x,y)$ (resp. $Z^{\sec}(x,y)$) when $f_{n}(k):=g_{n}(k)$ (resp. $f_{n}(k):=h_{n}(k)$). Then,
$$
\leqalignno{
Z^{\tan}(x,y)&=\sec(x+y)\cos(x-y);&(1.9)\cr
Z^{\sec}(x,y)&=\sec^2(x+y)\cos(x-y).&(1.10)\cr}
$$

As $Z^{\tan}(y,x)=Z^{\tan}(x,y)$ and 
$Z^{\sec}(y,x)=Z^{\sec}(x,y)$, this implies the following Corollary. 

\proclaim Corollary 1.3. The entries 
$g_{n}(k)$ and $h_{n}(k)$ have the symmetry property:
$$
g_{n}(k)=g_{n}(2n-k),\quad 
h_{n}(k)=h_{n}(2n-k)\quad (1\le k\le 2n-1).\leqno(1.11)
$$

In view of (1.7) and (1.8), two finite sets ${\goth A}_{2n-1}$ and ${\goth A}_{2n}$, of cardinalities $T_{2n-1}$ and $E_{2n}$, are to be found, together with a statistic, call it ``grn,'' defined on those sets with the property that
$$\leqalignno{
\sum_{\sigma\in {\goth A}_{2n-1}}x^{\grn\sigma}
&=\sum_{k}g_{n}(k)x^k;&(1.11)\cr}
$$

\goodbreak
$$\leqalignno{
\sum_{\sigma\in {\goth A}_{2n}}x^{\grn\sigma}
&=\sum_{k}h_{n}(k)x^k.&(1.12)\cr
}
$$
We shall use D\'esir\'e Andr\'e's old result [An1879, An1881], who introduced the notion of {\it alternating} permutation, as being a permutation 
$\sigma=\sigma(1)\sigma(2)\cdots\sigma(n)$ of $12\cdots n$ with the property that $\sigma(1)>\sigma(2)$,
$\sigma(2)<\sigma(3)$, $\sigma(3)>\sigma(4)$, etc. in an
alternating way. For each $n\ge 1$ let~${\goth A}_{n}$ denote the set of all alternating permutations of $12\cdots n$. He proved that $\#{\goth A}_{2n-1}=T_{2n-1}$, $\#{\goth A}_{2n}=E_{2n}$. 
The desired statistic ``grn'' is then the following.

\medskip
{\it Definition}.\quad
Let $\sigma=\sigma(1)\sigma(2)\cdots\sigma(n)$ be an alternating permutation from ${\goth A}_{n}$, 
so that $\sigma(i)=n$ for a certain~$i$ $(1\le i\le n)$. 
By convention, let $\sigma(0)=\sigma(n+1):=0$. Define the {\it {\bf gr}eater {\bf n}eighbor of~$n$ in~$\sigma$} to be 
$$\leqalignno{\noalign{\vskip-7pt}
{\rm grn}(\sigma)&:=\max\{\sigma(i-1), \sigma(i+1)\}.&(1.13)\cr
\noalign{\hbox{
Also, let
}}
{\goth A}_{n,k}&:=\{\sigma\in {\goth A}_{n}:\grn(\sigma)=k\}
\quad(0\le k\le n-1).&(1.14)\cr}
$$

\proclaim Theorem 1.4. Under the same assumptions as in Theorem~$1.1$ we have
$$\leqalignno{\qquad\quad
g_{n}(k)&=\#{\goth A}_{2n-1,k-1}\quad(n\ge 1,\, 1\le k\le 2n-1);&(1.15)\cr
h_{n}(k)&=\#{\goth A}_{2n,k}\quad(n\ge 1,\, 1\le k\le 2n-1).&(1.16)\cr}
$$

{\it Example}.\quad There are
$T_{3}=2$ alternating permutations of length~3, namely,  213 and 312, and ${\rm grn}(213)={\rm grn}(312)=1$, so that
$g_{2}(1)=\#{\goth A}_{3,0}=0$, $g_{2}(2)
=\#{\goth A}_{3,1}=2$,
$g_{2}(3)
=\#{\goth A}_{3,2}=0$; \quad there are 
$E_{4}=5$ alternating permutations of length~4, namely,
4132, 4231, 3142, 3241, 2143, and  ${\rm grn}(4132)=1$, ${\rm grn}(4231)=
{\rm grn}(3142)={\rm grn}(3241)=2$, ${\rm grn}(2143)=3$,
so that $h_{2}(1)=1$, $h_{2}(2)=3$, $h_{2}(3)=1$;
in accordance with the numerical values in Fig.~1.1.

\medskip
As $\#{\goth A}_{2n-1}=T_{2n-1}$, $\#{\goth A}_{2n}=E_{2n}$, following D\'esir\'e Andr\'e's result, it is now clear that
Theorem 1.1 is a consequence of Theorem 1.4. Thus, an analytical result is proved by combinatorial methods. 

In Proposition 2.1 it will be proved that 
$\#{\goth A}_{1,0} =1$,
$\#{\goth A}_{2n-1,0}=\#{\goth A}_{2n-1,2n-2}=0$ 
$(n\ge 2)$
and
$\#{\goth A}_{2n-1,1}=2\,T_{2n-3}$ $(n\ge 2)$;
also,
$\#{\goth A}_{2,1}=1$,
$\#{\goth A}_{2n,1}=\#{\goth A}_{2n,2n-1}=E_{2n-2}$ $(n\ge 2)$ and
$\#{\goth A}_{2n,2}=3\,E_{2n-2}$ $(n\ge 2)$.
In view of Theorem~1.4 this implies that conditions {\tt [tan2]} and {\tt [sec2]} are fulfilled as soon as 
conditions {\tt [tan1]} and {\tt [sec1]} hold, and then the following theorem.

\proclaim Theorem 1.5. The entries $(g_{n}(k))$ and $(h_{n}(k))$ given by $(1.15)$ and $(1.16)$
are also  solutions of the finite difference equation system $(1.2)$ when the initial values {\tt [tan2]} and {\tt [sec2]} are used, respectively.

\goodbreak
\smallskip
There are other combinatorial models which are also counted by tangent and secant numbers, or in a one-to-one correspondence with alternating permutations, in particular, the {\it labeled, binary, increasing, topological trees}, also called ``arbres binaires croissants complets" by Viennot [Vi88, chap.~3, p.~111]. The set of those trees having $n$ labeled nodes is denoted by ${\goth T}_{n}$. The statistic ``pom''
({\it {\bf p}arent {\bf o}f the {\bf m}aximum leaf\/}), introduced by Poupard [Po89] for her strictly ordered, binary trees can be extended to all of ${\goth T}_{n}$.
The usual bijection (see [Vi88]) $\gamma : {\goth T}_{n} \rightarrow {\goth A}_{n}$, called {\it projection}, has the property:
$\pom(t)=\grn(\gamma(t))$, as proved in Theorem~5.1.
We then have another combinatorial interpretation for the polynomials $\sum_{k}g_{n}(k)x^k$ and $\sum_{k}h_{n}(k)x^k$. 

\medskip
As such, the triangle $(g_{n}(k))$ $(n\ge 1,\,1\le k\le 2n-1)$ does not appear in Sloane's On-Line Encyclopedia of Integer Sequences [Sl06], but the triangle $(g_{n}(k)/2^{n-1})$ does under reference A008301 and is called Poupard's triangle, after her pioneering work on strictly ordered binary trees [Po89]. It is banal to verify that 
$2^{n-1}$ divides $g_{n}(k)$ when dealing with the combinatorial model  ${\goth T}_{n}$ for~$n$ odd.

In contrast to Christiane Poupard [Po89], who showed  that the
distribution of the strictly ordered binary trees satisfied the finite difference equation system $\Delta^2 f_{2n+1}(k)+2\,f_{2n-1}(k)=0$, we have used the multiplicative factor ``4"  in equation~(1.2) to make a
unified study of the tangent {\it and} secant cases and deal with
objects in one-to-one  correspondence with alternating permutations.
{\it Mutatis mutandis}, identities  (1.15), as  well
as (1.11) concerning the tangent numbers, are due to her.  
She obtains the generating function for her trees in the form:
$\sec((x+y)/\sqrt 2)\cos((x-y)/\sqrt 2)$ instead of (1.9).
However, the alternating permutation development in Sections~2 and~3,  identity (1.10)  and the combinatorial properties of the entries $h_{n}(k)$ are new.

The triangle of the $h_{n}(k)$'s appears in Sloane's [Sl06] as sequence A125053. It was deposited there by Paul D. Hanna. The entries have been calculated by using a procedure equivalent to (1.2) and the initial condition {\tt [tan2]}. No combinatorial interpretation is given and no generating function calculated.

\medskip
Theorem 1.4 will be proved in Section~3, once evaluations of some cardinalities such as $\#{\goth A}_{n,k}$ will be made, as done in Section~2. The proof of Theorem~1.2 is given in Section~4, together with further identities on the $g_{n}(k)$ and $h_{n}(k)$'s.

\bigskip
\centerline{\bf 2. Some special values}

\medskip
The evaluations of $\#{\goth A}_{2n-1,k-1}$ and  $\#{\goth A}_{2n,k}$ made in the next proposition for some values of $n$ and~$k$ will facilitate the derivation of the proof of Theorem~1.2. They also have their own combinatorial interests.

\proclaim Proposition 2.1. The following relations hold:
$$
\leqalignno{\noalign{\vskip-5pt}
\#{\goth A}_{1,0} =1,\quad 
\#{\goth A}_{2n-1,0}&=\#{\goth A}_{2n-1,2n-2}=0 
\quad(n\ge 2)
&(2.1)\cr
\#{\goth A}_{2n-1,1}=\#{\goth A}_{2n-1,2n-3}&=2\,T_{2n-3}\quad(n\ge 2);
&(2.2)\cr
\#{\goth A}_{2n-1,2}=\#{\goth A}_{2n-1,2n-4}&=4\,T_{2n-3}\quad(n\ge 3);&(2.3)\cr 
\#{\goth A}_{2n-1,3}=\#{\goth A}_{2n-1,2n-5}&=6\,T_{2n-3}
-8\,T_{2n-5}\quad(n\ge 3);&(2.4)\cr 
\sum_{k\ge 0}\#{\goth A}_{2n-1,k}&=T_{2n-1}
\quad(n\ge 1).&(2.5)\cr
\noalign{\vskip-6pt}
\#{\goth A}_{2,1}&=1;&(2.6)\cr
\#{\goth A}_{2n,1}=\#{\goth A}_{2n,2n-1}&=E_{2n-2}\quad (n\ge 2);&(2.7)\cr
\#{\goth A}_{2n,2}=\#{\goth A}_{2n,2n-2}&=3\,E_{2n-2}\quad (n\ge 2);&(2.8)\cr
\#{\goth A}_{2n,3}=\#{\goth A}_{2n,2n-3}&=5\,E_{2n-2}-4\,E_{2n-4}\quad (n\ge 2);&(2.9)\cr
\sum_{k\ge 1}\#{\goth A}_{2n,k}&=E_{2n}
\quad(n\ge 1).&(2.10)\cr}
$$

\goodbreak
{\it Proof}.\quad
(2.1)\quad The set ${\goth A}_{1}$ is the singleton~1 and ${\rm grn}(1)=0$ by definition, so that $\#{\goth A}_{1,0}=1$. For $n\ge 2$
all alternating permutations from ${\goth A}_{2n-1}$ have a ``grn'' at least equal to~1. Hence, $\#{\goth A}_{2n-1,0}=0$. Finally, each alternating permutation of length $(2n-1)$ $(n\ge 2)$ contains neither the factor $(2n-2)(2n-1)$, nor $(2n-1)(2n-2)$. Hence,
$\#{\goth A}_{2n-1,2n-2}=0$.

(2.2)\quad When $n\ge 2$, each alternating permutation from ${\goth A}_{2n-1,1}$ starts with  $(2n-1)\, 1$,  or ends with $1\,(2n-1)$. After removal of those two letters, there remains an alternating permutation on $\{2,3,\ldots,2n-2\}$. Hence, $\#{\goth A}_{2n-1,1}=2\,T_{2n-3}$. Next, each permutation from ${\goth A}_{2n-1,2n-3}$ 
must contain, either the three-letter factor $(2n-1)(2n-3)(2n-2)$, or $(2n-2)(2n-3)(2n-1)$. The removal of the factor 
$(2n-1)(2n-3)$ (resp. $(2n-3)(2n-1)$) yields an alternating permutation of the set $\{1,2,\ldots,(2n-4),(2n-2)\}$, of cardinality $(2n-3)$. This proves relation
$\#{\goth A}_{2n-1,2n-3}=2\,T_{2n-3}$. 

(2.3)\quad Start with an alternating permutation on
$\{1,3,4,\ldots,2n-2\}$, then having $(2n-3)$ elements. There are four possibilities to generate a permutation from ${\goth A}_{2n-1,2}$: (1) insert $(2n-1)\,2$ to the left; (2) insert $2\,(2n-1)$ to the right; (3) insert $2 (2n-1)$ just before~1; (4) insert $(2n-1)\ 2$ just after~1. For the second identity in (2.3) proceed in the same way: in each alternating permutation on $\{1,2,\ldots, 2n-1\}\setminus\{2n-4,2n-1\}$  the two letters $(2n-3)$, $(2n-2)$ are necessarily local maxima. There are four possibilities to obtain a permutation from ${\goth A}_{2n-1,2n-4}$: insert $(2n-1)\,(2n-4)$ just before, either $(2n-3)$, or $(2n-2)$; also insert $(2n-4)\,(2n-1)$ just after, either $(2n-3)$, or $(2n-2)$.

(2.4)\quad Each permutation from ${\goth A}_{2n-1,3}$ containing the factor $2\,1$ (resp. $1\,2$) starts with $2\,1$ (resp. ends with $1\,2$). Dropping the factor $2\,1$ (resp. $1\,2$) and subtracting~2 from the remaining letters yields an alternating permutation from ${\goth A}_{2n-3,1}$.
There are then $2(2\,T_{2n-5})$ permutations from ${\goth A}_{2n-1,3}$ containing, either $2\,1$, or $1\,2$.

If a permutation from ${\goth A}_{2n-1,3}$ contains neither one of those two factors, it has one of the {\it six} properties: it starts with $(2n-1)\,3$, or contains one of the three-letter factor $1\,(2n-1)\,3$, $3\,(2n-1)\,1$, $2\,(2n-1)\,3$, $3\,(2n-1)\,2$, or still ends with $3\,(2n-1)$. After removal of the two-letter factor $(2n-1)\,3$ or $3\,(2n-1)$ there remains an alternating permutation on $\{1,2,4,\ldots,(2n-2)\}$ {\it not starting with} $2\,1$
and {\it not ending with} $1\,2$. There are then 
$6(T_{2n-3}-2\,T_{2n-5})$ such permutations. Altogether,
$\#{\goth A}_{2n-1,3}=4\,T_{2n-5}+6(T_{2n-3}-2\,T_{2n-5})=
6\,T_{2n-3}-8\,T_{2n-5}$.

The proof of the second identity in (2.4) follows a different pattern. If the letter $(2n-4)$ is a local minimum (i.e., less than its two adjacent letters) in a permutation~$\sigma$ from ${\goth A}_{2n-1,2n-5}$, then~$\sigma$ necessarily contains one of the four five-letter factors $(2n-1)(2n-5)(2n-2)(2n-4)(2n-3)$,
$(2n-1)(2n-5)(2n-3)(2n-4)(2n-2)$,\ 
$(2n-2)(2n-4)(2n-3)(2n-5)(2n-1)$,
$(2n-3)(2n-4)(2n-2)(2n-5)(2n-1)$. 
Replacing this five-letter factor by $(2n-5)$ yields a permutation from ${\goth A}_{2n-5}$. Thus, there are $4\,T_{2n-5}$ permutations 
from ${\goth A}_{2n-1,2n-5}$ in which $(2n-4)$ is a local minimum.

In the other permutations from ${\goth A}_{2n-1,2n-5}$ all the four letters $(2n-4)$, $(2n-3)$, $(2n-2)$, $(2n-1)$ are local maxima (i.e., greater than their adjacent letters). Let 
${\goth A}_{2n-1,2n-5}'$ be the set of those permutations. When the two-letter factor $(2n-1)(2n-5)$ or $(2n-5)(2n-1)$ is deleted from such a permutation, there remains a permutation on
$\{1,2, \ldots,(2n-1)\}\setminus\{(2n-5), (2n-1)\}$ in which
the third largest letter $(2n-4)$ is not a local minimum. Let
${\goth A}_{2n-3}''$ be the set of those permutations.
But the alternating permutations on the latter set in which $(2n-4)$ is a local minimum necessarily contain the three-letter factor $(2n-3)(2n-4)(2n-2)$ or $(2n-2)(2n-4)(2n-3)$. There are then $2\,T_{{2n-5}}$ such permutations. Hence, $\#{\goth A}_{2n-3}''=T_{2n-3}-2\,T_{2n-5}$.
To obtain a permutation from ${\goth A}_{2n-1,2n-5}'$ it suffices to start from a permutation~$\sigma''$ from 
${\goth A}_{2n-3}''$ and insert $(2n-1)(2n-5)$ (resp.
$(2n-5)(2n-1)$) just before (resp. just after) each one of the three letters $(2n-4)$, $(2n-3)$, $(2n-2)$ (which are all local maxima). There are then $6(T_{2n-3}-2\,T_{2n-5})$ such permutations. Altogether, $\#{\goth A}_{2n-1,2n-5}
=4\,T_{2n-5}+6(T_{2n-3}-2\,T_{2n-5})=6\,T_{2n-3}-8\,T_{2n-5}$.
No comment for (2.5) and (2.6).

(2.7)\quad Simply note that the only alternating permutations  from ${\goth A}_{2n,1}$ and ${\goth A}_{2n,2n-1}$ are, respectively, of the form:\quad $(2n)\, 1\, \sigma(3)\cdots \sigma(2n)$ and \hfil\break $\sigma(1)\sigma(2)\cdots (2n)\,(2n-1)$.

(2.8)\quad Same proof as for (2.3): start with an alternating permutation on $\{1,3,4,\ldots,(2n-1)\}$. There are exactly three possibilities to generate a permutation from ${\goth A}_{2n,2}$:
insert $(2n)\,2$ to the left, or just after the letter~1, or still insert $2\,(2n)$ just before the letter~1.
For the second identity start with a permutation on $\{1,2,\ldots, 2n\}\setminus\{2n-2,2n\}$ and insert $(2n)\,(2n-2)$ either to the right, or just before $(2n-1)$, or still insert $(2n-2)\,(2n)$ just after $(2n-1)$.

(2.9)\quad Each permutation from ${\goth A}_{2n,3}$ containing the factor $2\,1$ is necessarily of the form $\sigma=2\,1\,\sigma(3)\cdots\sigma(2n)$, so that the alternating permutation $\sigma':=(\sigma(3)-2)\cdots(\sigma(2n)-2)$ belongs to ${\goth A}_{2n-2,1}$.
There are then $E_{2n-4}$ such permutations. If a permutation from
${\goth A}_{2n,3}$ does not contain $2\,1$, it has one of the {\it five} properties: it starts with $(2n)\,3$, or contains one of the three-letter factor $1\,(2n)\,3$, $3\,(2n)\,1$, $2\,(2n)\,3$, $3\,(2n)\,2$. After removal of the two-letter factor $(2n)\,3$ or $3\,(2n)$ there remains an alternating permutation on $\{1,2,4,\ldots,(2n-1)\}$ {\it not starting with} $2\,1$. There are then 
$5(E_{2n-2}-E_{2n-4})$ such permutations. Altogether,
$\#{\goth A}_{2n,3}=E_{2n-4}+5(E_{2n-2}-E_{2n-4})=
5E_{2n-2}-4E_{2n-4}$.

The proof for the second identity in (2.9) is quite similar. Each permutation from ${\goth A}_{2n,2n-3}$ containing the factor $(2n-1)(2n-2)$ necessarily ends with the four-letter factor $(2n)(2n-3)(2n-1)(2n-2)$. There are then $E_{2n-4}$ such permutations. The other permutations from ${\goth A}_{2n,2n-3}$ contain one of the four three-letter factors $(2n)(2n-3)(2n-2)$, $(2n-2)(2n-3)(2n)$, $(2n)(2n-3)(2n-1)$, $(2n-1)(2n-3)(2n)$, or ends with $(2n)(2n-3)$. After removal of the two-letter factor $(2n)(2n-3)$ or $(2n-3)(2n)$ there remains an alternating permutation on $\{1,2,\ldots,(2n-4),(2n-2),(2n-1)\}$, not ending with the two-letter factor $(2n-1)(2n-2)$. There are $E_{2n-2}-E_{2n-4}$ such permutations. Altogether, $\#{\goth A}_{2n,2n-3}=E_{2n-4}+5(E_{2n-2}-E_{2n-4})$. 

No comment for (2.10).\qed

\bigskip
\centerline{\bf 3. Proof of Theorem 1.4}
\medskip
Let $a_{n}(k):=\#{\goth A}_{2n-1,k-1}$ and
$b_{n}(k):=\#{\goth A}_{2n,k}$.
From Proposition~2.1 it follows that the initial conditions {\tt [tan1]} and {\tt [tan2]} hold when $f_{n}(k)=a_{n}(k)$, and {\tt [sec1]} and also {\tt [sec2]} when $f_{n}(k)=b_{n}(k)$. It remains to prove that in each case (1.2) holds.

By means of identities (2.2)--(2.4) and (2.7)--(2.9)
we easily verify that (1.2) holds for both
$a_{n}(k)$ and $b_{n}(k)$ when $n=2,3$ and $1\le k\le 2n-3$. It also holds for $a_{n}(k)$ when $n\ge 4$ and $k=1,2,2n-4,2n-3$,
and for $b_{n}(k)$ when $n\ge 4$ and $k=1,2n-3$.

What is left to prove is:
$\Delta^2a_{n}(k)+4a_{n-1}(k)=0$, that is, 
$\Delta^2{\goth A}_{2n-1,k-1}+4\,{\goth A}_{2n-3,k-1}=0$  
for $n\ge 4$ and $3\le k\le 2n-5$---by identifying
each finite set with its cardinality---and also
$\Delta^2b_{n}(k)+4b_{n-1}(k)=0$, that is, 
$\Delta^2{\goth A}_{2n,k}+4\,{\goth A}_{2n-2,k}=0$  
for $n\ge 4$ and $2\le k\le 2n-4$; altogether,
$$
\Delta^2{\goth A}_{n,k}+4\,{\goth A}_{n-2,k}=0
\quad {\rm for}\ n\ge 7\ {\rm and}
\vtop{\hbox{$\ 2\le k \le n-4\ (n\ {\rm even})$}
\hbox{$\ 2\le k \le n-5\  (n\ {\rm odd}).$}}
\leqno(2.11)
$$

Let $v=y_{1}\cdots y_{m}$ be a nonempty word with {\it distinct} letters from the set $\{0,1,2,\ldots,n\}$ and 
$\widetilde v=y_{m}\cdots y_{1}$ be its mirror-image. If $m=1$ and 
$y_{1}=0$, let $[v]=[0]$ be the empty set.
If $m\ge 2$ and $y_{1}=0$ (resp. $y_{m}=0$), let $[v]$ be the set of all alternating permutations from~${\goth A}_{n}$, if any, whose left factors are equal to $y_{2}\cdots y_{m}$, or whose right factors are equal to $y_{m}\cdots y_{2}$. When $y_{1}\ge 1$, let $[v]$ be the set of all alternating permutations from~${\goth A}_{n}$, if any,
containing, either the factor $v$, or the factor $\widetilde v$. Finally, let $[\,\widetilde v\,]:=[v]$.

Using those notations we get
$$\leqalignno{
{\goth A}_{n,k}&=\sum_{0\le y\le k-1}[ynk]\cr
&=\sum_{0\le y\le k-1}[ynk(k+1)]
+\sum_{\scriptstyle 0\le y\le k-1
\atop\scriptstyle k+2\le z\le n-1}[ynkz]
+\sum_{1\le y\le k-1}[ynk0];\cr
\noalign{\vskip-5pt}
{\goth A}_{n,k+1}&=\sum_{0\le y\le k}[yn(k+1)]\cr
&=[kn(k+1)]+\sum_{\scriptstyle 0\le y\le k-1
\atop\scriptstyle k+2\le z\le n-1}[yn(k+1)z]
+\sum_{1\le y\le k-1}[yn(k+1)0].
\cr
}
$$
The transposition $(k,k+1)$ maps the set $[ynkz]$ onto the set $[yn(k+1)z]$ for $z\in \{k+2,\ldots,n-1\} \cup \{0\}$, so that we may write
$$\leqalignno{\Delta\,{\goth A}_{n,k}=
{\goth A}_{n,k+1}-{\goth A}_{n,k}
&=[kn(k+1)]-\sum_{0\le y\le k-1}[ynk(k+1)]\cr
&=[kn(k+1)]-\sum_{\scriptstyle 0\le y_{1},y_{2}\le k-1
\atop\scriptstyle y_{1}\not=y_{2}}[y_{1}nk(k+1)y_{2}];\cr
}$$
$$
\displaylines{\noalign{\vskip-18pt}
\quad\Delta\,{\goth A}_{n,k+1}=
{\goth A}_{n,k+2}-{\goth A}_{n,k+1}\hfill\cr
\qquad{}=[(k+1)n(k+2)]-\sum_{0\le y\le k}[yn(k+1)(k+2)]\hfill\cr
\qquad{}=[(k+1)n(k+2)]-[kn(k+1)(k+2)]\hfill\cr
\kern1cm{}-
\sum_{0\le y\le k-1}[yn(k+1)(k+2)k]
-\kern-10pt \sum_{\scriptstyle 0\le y_{1},y_{2}\le k-1
\atop\scriptstyle y_{1}\not=y_{2}}[y_{1}n(k+1)(k+2)y_{2}]
.\hfill\cr
}
$$
For $2\le k\le n-4$ the permutation ${\;\;k\quad k+1\;\;k+2\choose \!k+1\;\; k+2\;\;k}$ maps $[y_{1}nk(k+1)y_{2}]$
onto $[y_{1}n(k+1)(k+2)y_{2}]$ in a bijective manner. Hence,
$$
\leqalignno{
\Delta^2{\goth A}_{n,k}&=\Delta\,{\goth A}_{n,k+1}
-\Delta\,{\goth A}_{n,k}\cr
&=[(k+1)n(k+2)]-[kn(k+1)(k+2)]\cr
&\kern1cm {}-\sum_{0\le y\le k-1}[yn(k+1)(k+2)k]
-[kn(k+1)].\cr}
$$
But
$$\leqalignno{
[kn(k+1)]&=[(k+2)kn(k+1)]+[kn(k+1)(k+2)]\cr
&\kern2.5cm{}+\sum_{\scriptstyle z_{1},z_{2}\in
\{k+3,\ldots,n-1\}\cup\{0\}\atop
\scriptstyle z_{1}\not=z_{2}}
[z_{1}kn(k+1)z_{2}].\cr
}
$$
Again, the permutation ${\;\;k\quad k+1\;\;k+2\choose \!k+1\;\; k+2\;\;k}$ maps the last sum onto the set $[(k+1)n(k+2)]$. Altogether, as
$[(k+2)kn(k+1)]=[kn(k+1)(k+2)]$, we have
$$
\Delta^2{\goth A}_{n,k}=
-3\;[kn(k+1)(k+2)]-\sum_{0\le y\le k-1}[yn(k+1)(k+2)k].
$$

When removing the factor $(k+1)(k+2)$ and replacing each integer $z\geq k+2$
by $(z-2)$, in each alternating permutation, both sets
$[kn(k+1)(k+2)]$ and
$\sum\limits_{0\le y\le k-1}[yn(k+1)(k+2)k]$
are transformed into ${\goth A}_{n-2,k}$;
so that 
$\Delta^2{\goth A}_{n,k}=-4\,{\goth A}_{n-2,k}$.

\medskip
\centerline{\bf 4. The bivariate generating functions}

\medskip
Let $f=(f_n(k))$ $(n\ge 1,\,1\le k\le 2n-1)$ be the family of rational
numbers, as displayed in (1.1), that satisfies the
finite-difference equation system~(1.2) under the initial conditions
{\tt [tan2]} of {\tt [sec2]}. We know that the system has then a unique solution. With the triangle~$f$ associate the infinite  matrix
$$
\Gamma\!=\!(\gamma_{ij})_{(i\ge 0,j\ge 0)}
:=\petitematrice{f_1(1)&0&f_2(3)&0&f_3(5)&0&f_4(7)&\cdots\cr
0&f_2(2)&0&f_3(4)&0&f_4(6)&\cdots\cr
f_2(1)&0&f_3(3)&0&f_4(5)&\cdots\cr
0&f_3(2)&0&f_4(4)&\cdots\cr
f_3(1)&0&f_4(3)&\cdots\cr
0&f_4(2)&\cdots\cr
f_4(1)&\cdots\cr}.\leqno(4.1)
$$
In other words, define $\gamma_{ij}:=0$ when $i+j$ is odd, and
$\gamma_{ij}:=f_{n}(k)$ with $k:=j+1$, $2n=2+i+j$ when $i+j$ is even. For $i+j$ even the mapping $(i,j)\mapsto (n,k)$ is one-to-one, the reverse mapping being for $n\ge 1$, $1\le k\le 2n-1$ given by
$i=2n-1-k$, $j=k-1$.

\goodbreak
In terms of the entries $\gamma_{ij}$ relation (1.2) may be
written in the form
$$\leqalignno{
\gamma_{i,j}&=2\,\gamma_{i-1,j-1}+{1\over
2}(\gamma_{i-1,j+1} +\gamma_{i+1,j-1})\quad (i\ge 1,\,j\ge
1);&(4.2)\cr
\gamma_{ij}&=0, \qquad {\rm if}\ i+j\ {\rm odd}.&(4.3)\cr
}
$$
Furthermore, the full matrix
$\Gamma=(\gamma_{i,j})$
$(i\ge 0,\,j\ge 0)$ is completely determined as soon as its
{\it first row} $(\gamma_{0,j})$ $(j\ge 0)$ and {\it first column}
$(\gamma_{i,0})$ $(i\ge 0)$ are known. Let $f\mapsto \Gamma$
denote the above correspondence between those  triangles  and
matrices. 

Let
$Z(x,y):=\displaystyle\sum_{i\ge 0,\,j\ge
0}\gamma_{i,j}{x^i\over i!}{y^j\over j!}$.
It is easily verified that $Z(x,y)$ satisfies the partial 
differential equation
$${\partial^2Z(x,y)\over \partial x\,\partial y}
=2\,Z(x,y)+{1\over 2}{\partial^2Z(x,y)\over\partial x^2}+
{1\over 2}{\partial^2Z(x,y)\over\partial y^2},\leqno(4.4)
$$
if and only if the
coefficients~$\gamma_{i,j}$ satisfy relation (4.2). Hence, 
$Z(x,y)$ is fully determined by (4.4) and by  the generating
functions
$Z(x,0)=\sum\limits_{i\ge 0}\gamma_{i,0}\,x^i/i!$ and
$Z(0,y)=\sum\limits_{j\ge 0}\gamma_{0,j}\,y^j/j!$ for the first
column and first row  of the matrix~$\Gamma$.

But for any given formal power series in one variable
$f(x)=1\!+\!\sum\limits_{n\ge
1}f_{2n}\,\displaystyle{x^{2n}\over(2n)!}$ it can be also verified
that the bivariate formal power series
$$\
Z(x,y)=\sum_{i\ge 0,\,j\ge
0}\gamma_{i,j}{x^i\over i!}{y^j\over j!}
=f(x+y)\sec(x+y)\cos(x-y)\qquad\leqno(4.5)
$$
satisfies (4.4) and that the generating functions for its first
column and first row are given by $f(x)$ and $f(y)$,
respectively. This proves the following proposition.

\proclaim Proposition 4.1. Let $f(x)=1\!+\!\sum\limits_{n\ge
1}f_{2n}\,\displaystyle{x^{2n}\over(2n)!}$ be given and
$\Gamma=(\gamma_{ij})$ $(i\ge 0,\,j\ge 0)$ be an infinite
matrix,  whose entries  satisfy relations $(4.2)$ and $(4.3)$, on
the one hand, and such that $\gamma_{0,0}=1$,
$\gamma_{2n+1,0}=\gamma_{0,2n+1}=0$ for $n\ge 0$
and  $\gamma_{2n,0}=\gamma_{0,2n}=f_{2n}$ for $n\ge 1$, on the
other hand. Then, identity  $(4.5)$ holds. 

Using the correspondence $\gamma_{ij}\leftrightarrow f_{n}(k)$ above mentioned, the series
$Z(x,y)$ can be rewritten
$$
Z(x,y)=1+\sum_{n\ge 1}\sum_{1\le k\le 2n+1}
f_{n+1}(k){x^{2n+1-k}\over (2n+1-k)!}
{y^{k-1}\over (k-1)!},\leqno(4.6)
$$
which is then equal to $f(x+y)\sec(x+y)\cos(x-y)$ under the assumptions of the previous proposition.

\goodbreak
Now, consider the two  triangles  described  in Fig.~1.1 and 
let $\Gamma^{\tan}=(\gamma_{ij}^{\tan})$  and
$\Gamma^{\sec}=(\gamma_{ij}^{\sec})$ be the two 
$\Gamma$-matrices attached to them:
$$\displaylines{
\Gamma^{\tan}\!=\!(\gamma_{ij}^{\tan})=
\petitematrice{g_{1}(1)&0&g_{2}(3)&0&g_{3}(5)&0&g_{4}(7)&\cdots\cr
0&g_2(2)&0&g_3(4)&0&g_4(6)&\cdots\cr
g_{2}(1)&0&g_3(3)&0&g_4(5)&\cdots\cr
0&g_3(2)&0&g_4(4)&\cdots\cr
g_{3}(1)&0&g_4(3)&\cdots\cr
0&g_4(2)&\cdots\cr
g_{4}(1)&\cdots\cr};\cr
\Gamma^{\sec}\!=\!(\gamma_{ij}^{\sec})=
\petitematrice{h_{1}(1)&0&h_{2}(3)&0&h_{3}(5)&0&h_{4}(7)&\cdots\cr
0&h_2(2)&0&h_3(4)&0&h_4(6)&\cdots\cr
h_{2}(1)&0&h_3(3)&0&h_4(5)&\cdots\cr
0&h_3(2)&0&h_4(4)&\cdots\cr
h_{3}(1)&0&h_4(3)&\cdots\cr
0&h_4(2)&\cdots\cr
h_{4}(1)&\cdots\cr}.\cr
}
$$
The exponential generating
function for the first row and first column of
$\Gamma^{\tan}$ is equal to
$f(x)=1$. On the other hand, as $h_{n}(1)=
h_{n}(2n+1)=\#{\goth A}_{2n,1}=E_{2n-2}$ 
by {\tt [sec2]}, (1.8) and (1.16),
the exponential generating
function for the first row and first column of
$\Gamma^{\sec}$ is equal to
$h_{1}(1)+h_{2}(1)x^2/2!+h_{3}(1)x^4/4!+\cdots
=E_{0}+E_{2}x/2!+E_{4}x^4/4+\cdots=\sec(x)$.
Theorem~1.2 is then a
consequence of the previous Proposition.

\medskip

When $(x,y)$ is equal to $(x,x)$, then to $(x,-x)$ in (1.9)
and  (1.10), we obtain: 
$Z^{\tan}(x,x)=\sec(2x)$; $Z^{\tan}(x,-x)=\cos(2x)$;
$Z^{\sec}(x,x)=\sec^2(2x)=1+\sum_{n\ge 1}4^nT_{2n+1}x^{2n}/(2n)!$
and $Z^{\sec}(x,-x)=\cos(2x)$. Looking for the coefficients of
$x^{2n}/(2n)!$ on both sides in the first
(resp. last) two formulas yields four further identities

$$
\leqalignno{\noalign{\vskip-16pt}
\sum_{1\le k\le  2n+1}{2n\choose k-1}g_{n+1}(k)&=4^{n}E_{2n}\quad 
(n\ge1);&(4.7)\cr
\sum_{1\le k\le  2n+1}(-1)^k{2n\choose
k-1}g_{n+1}(k)&=(-1)^n4^{n}\quad 
(n\ge1);&(4.8)\cr
\sum_{1\le k\le  2n+1}{2n\choose
k-1}h_{n+1}(k)&=4^{n}T_{2n+1}\quad  (n\ge1)&(4.9)\cr
\sum_{1\le k\le  2n+1}(-1)^k{2n\choose
k-1}h_{n+1}(k)&=(-1)^n4^{n}\quad 
(n\ge1).&(4.10)\cr}
$$
The last two ones are mentioned in Sloane's Encyclopedia [Sl07]
(sequence A125053) without proofs.

\newbox\boxarbre
\medskip
\setbox\boxarbre=\vbox{\vskip
20mm\offinterlineskip 
\centerput(33,20){$7$}
\centerput(10,21){$5$}
\centerput(20,21){$8$}
\centerput(15,12){$4$}
\centerput(42,16){$3$}
\centerput(0,11){$6$}
\centerput(30,7){$2$}\centerput(45,10){$\mapsto$}
\centerput(5,2){$1$}
\segment(0,10)\dir(1,-1)\long{5}
\segment(5,5)\dir(5,1)\long{25}
\segment(30,10)\dir(2,1)\long{10}
\segment(40,15)\dir(-2,1)\long{7}
\segment(30,10)\dir(-3,1)\long{15}
\segment(15,15)\dir(1,1)\long{5}
\segment(15,15)\dir(-1,1)\long{5}
}

\newbox\boxarbreb
\medskip
\setbox\boxarbreb=\vbox{\vskip
20mm\offinterlineskip 
\centerput(33,20){$8$}
\centerput(10,21){$4$}
\centerput(20,21){$7$}
\centerput(15,12){$3$}
\centerput(42,16){$6$}
\centerput(0,11){$5$}
\centerput(30,7){$2$}
\centerput(5,2){$1$}
\segment(0,10)\dir(1,-1)\long{5}
\segment(5,5)\dir(5,1)\long{25}
\segment(30,10)\dir(2,1)\long{10}
\segment(40,15)\dir(-2,1)\long{7}
\segment(30,10)\dir(-3,1)\long{15}
\segment(15,15)\dir(1,1)\long{5}
\segment(15,15)\dir(-1,1)\long{5}
}

\bigskip
\centerline{\bf 5. Alternating permutations and binary trees}

\medskip

In this Section the traditional vocabulary on trees, such
as node, leaf, child, root, \dots\ is used. In particular, when
a node is not a leaf, it is said to be an {\it internal node}.

\smallskip
{\it Definition}.\quad
An {\it $n$-labeled, binary, increasing, topological tree}
is defined by the following axioms:

(1) it is a {\it labeled\/} tree with $n$ nodes, labeled
$1,2,\ldots, n$; the node labeled~1 is called the {\it root};

(2) each node has no child (then called a {\it leaf\/}), or
one child, or two children;

(3) the label of each node is smaller than the label of its children, if any;

(4) the tree is planar and each child of a node is, either on the
left (it is then called the {\it left child}), or on the right (the
{\it right child}); 

(5) when $n$ is odd, each node is, either a leaf, or
a node with two children; when $n$ is even, each node is, either a leaf, or a node with two children, except the rightmost node
(uniquely defined) which has one {\it left} child, but no right child.
It will be referred to as being the {\it one-son child}.

\smallskip
Each such binary tree $t$ may be drawn on a Euclidean plane:
the root has coordinates $(0,0)$, the left son of the root
$(-1,1)$ and the right son
$(1,1)$, the grandsons $(-3/2,2)$, $(-1/2,2)$,
$(1/2,2)$, $(3/2,2)$, respectively, the great-grandsons
$(-7/4 ,3)$, $(-5/4,3)$, \dots~, $(7/4,3)$, etc. With this
convention all the nodes have different abscissas. Let~$t$ have
$n$ nodes and make the orthogonal projections of those nodes
on a horizontal axis. Writing the labels of the projected
$n$ nodes yields a
permutation~$\sigma=\sigma(1)\sigma(2)\cdots\sigma(n)$ of
$1\,2\,\cdots\,n$. We say that $\sigma$ is the {\it projection} of~$t$ and~$t$ the {\it spreading out} of~$\sigma$. Moreover, 
$\sigma$ is {\it alternating}. For instance, the two trees $t_1$ and $t_2$
in Fig.~1.2 are labeled, binary, increasing,
topological trees, with~7 and~8 nodes, respectively. Their
projections $\sigma_1=6\,1\,5\,4\,7\,2\,3$
and~$\sigma_2=6\,1\,5\,4\,8\,2\,7\,3$ are
alternating. 

For each $n\ge 1$ let ${\goth T}_n$  be the set
of all $n$-labeled, binary, increasing, topological trees. Then,
$t\mapsto\sigma$ is a bijection of ${\goth T}_n$ onto ${\goth A}_n$,
so that we also have:
$\#{\goth A}_{2n+1}=T_{2n+1}$,
$\#{\goth A}_{2n}=E_{2n}$. Each tree~$t$ from~${\goth
T}_n$ is said to be {\it tangent} (resp. {\it secant}),
if~$n$ is odd (resp. even).

\newbox\boxarbre
\newbox\boxarbrea

\setbox\boxarbre=\vbox{\vskip
25mm\offinterlineskip 
\centerput(10,21){$5$}
\centerput(20,21){$7$}
\centerput(15,12){$4$}
\centerput(30,16){$3$}
\centerput(0,11){$6$}
\centerput(25,7){$2$}
\centerput(5,2){$1$}
\segment(0,10)\dir(1,-1)\long{5}
\segment(5,5)\dir(4,1)\long{20}
\segment(25,10)\dir(1,1)\long{5}
\segment(25,10)\dir(-2,1)\long{10}
\segment(15,15)\dir(1,1)\long{5}
\segment(15,15)\dir(-1,1)\long{5}
\segment(-5,0)\dir(1,0)\long{40}
\centerput(0,-4){$6$}
\centerput(5,-4){$1$}
\centerput(10,-4){$5$}
\centerput(15,-4){$4$}
\centerput(20,-4){$7$}
\centerput(25,-4){$2$}
\centerput(30,-4){$3$}
\centerput(-5,-4){$\sigma_1={}$}
\leftput(-10,15){$t_1={}$}
}

\setbox\boxarbrea=\vbox{\vskip
25mm\offinterlineskip 
\centerput(10,21){$5$}
\centerput(20,21){$8$}
\centerput(15,12){$4$}
\centerput(30,21){$7$}
\centerput(0,11){$6$}
\centerput(25,7){$2$}
\centerput(5,2){$1$}
\centerput(35,12){$3$}
\segment(0,10)\dir(1,-1)\long{5}
\segment(5,5)\dir(4,1)\long{20}
\segment(25,10)\dir(2,1)\long{10}
\segment(25,10)\dir(-2,1)\long{10}
\segment(15,15)\dir(1,1)\long{5}
\segment(15,15)\dir(-1,1)\long{5}
\segment(35,15)\dir(-1,1)\long{5}
\segment(-5,0)\dir(1,0)\long{45}
\centerput(0,-4){$6$}
\centerput(5,-4){$1$}
\centerput(10,-4){$5$}
\centerput(15,-4){$4$}
\centerput(20,-4){$8$}
\centerput(25,-4){$2$}
\centerput(30,-4){$7$}
\centerput(35,-4){$3$}
\centerput(-5,-4){$\sigma_2={}$}
\leftput(-10,15){$t_2={}$}
}

\vskip-10pt
$$
\box\boxarbre\hskip50mm\box\boxarbrea\hskip3.3cm
$$

\medskip
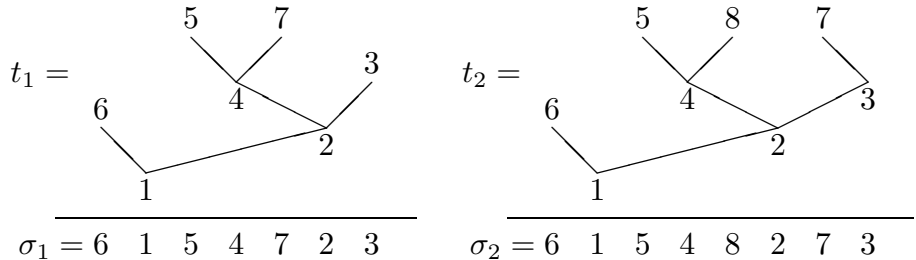
\centerline{Fig. 1.2. Tangent, secant trees and alternating
permutations.}

\medskip
The two statistics ``emc'' (``{\bf e}nd {\bf o}f {m}inimal {\bf c}hain'') and ``pom'' (``{\bf p}arent {\bf o}f {\bf m}aximum {l}eaf'') we now define
on each set ${\goth T}_n$ have been introduced by Christiane
Poupard [Po89] in her study of the strictly ordered binary trees and provide two other combinatorial interpretations for the entries $g_{n}(k)$ and $h_{n}(k)$. Their definitions are also valid for all binary increasing trees, in particular, for secant and tangent trees. Let $n\ge 2$ and~$t$ be a binary increasing tree, with~$n$ nodes labeled $1,2,\ldots,n$. Let~$a$ be the label of  an internal node.
If the node has two children labeled $b$ and $c$,
define $\min a:=\min\{b,c\}$; if it has one child~$b$, let $\min a:=b$.
The {\it minimal chain} of~$t$ is defined to~be~the sequence
$a_1\rightarrow a_2\rightarrow a_3\rightarrow\cdots
\rightarrow a_{j-1}\rightarrow a_j$, with the following properties:

(i) $a_1=1$ is the label of the root;

(ii) for
$i=1,2,\ldots,j-1$ the $i$-th term~$a_{i}$ is the label of an
internal node and $a_{i+1}=\min a_i$;

(iii) $a_j$ is the label of a leaf. 

\noindent 
Define 
the {\it end of the  minimal  chain\/} in~$t$
to be:
$\eoc(t):=a_j$.
As~$t$ is increasing, there is a  unique leaf with label~$n$. If that
leaf is incident to a node labeled~$k$, define the ({\it parent of the maximum leaf\/}) in~$t$ to be:
$\pom(t):=k$. By convention, $\eoc(t)=1$ and $\pom(t)=0$ for the unique~$t\in
{\goth T}_1$.

The minimal chain of the tree $t_{1}$ (resp. $t_{2}$) displayed in Fig.~1.2 is $1\rightarrow 2\rightarrow
3$ (resp. $1\rightarrow 2\rightarrow
3\rightarrow 7$). Then $\eoc(t_{1})=3$, $\eoc(t_{2})=7$. Also, $\pom(t_{1})=4$ and $\pom(t_{2})=4$.

\proclaim Theorem 5.1. Let $\sigma=\sigma(1)\sigma(2)\cdots \sigma(n)$ be the projection of the $n$-labeled, binary, increasing, topological tree~$t$. Then $\pom(t)={\rm grn}(\sigma)$. In other words, the parent of the maxium leaf in~$t$ is the greater neighbor of~$n$ in~$\sigma$.

\goodbreak
{\it Proof}.\quad
Let $\sigma(i)=n$ with $2\le i\le n-1$. The parent of the node labeled~$n$ in~$t$ is, either the node labeled~$\sigma(i-1)$, or the node labeled $\sigma(i+1)$. Let $\sigma(j)$ be the label of the node of the common ancester of the previous two nodes in~$t$. Then, $\sigma(j)\le \min \{\sigma(i-1),\,\sigma(i+1)\}$, $j\not=i$  and $i-1\le j\le i+1$. Hence, either $\sigma(j)=\sigma(i-1)<\sigma(i+1)$, or $\sigma(j)=\sigma(i+1)<\sigma(i-1)$.
In the first (resp. the second) case the parent of~$n$ is $\sigma(i+1)$ (resp. $\sigma(i-1)$) and ${\rm grn}(\sigma)=\max\{\sigma(i-1),\sigma(i+1)\}$.\qed

\medskip
In [Ha12] an explicit bijection of ${\goth T}_n$ onto itself is constructed that maps the statistic ``$\eoc-1$'' onto the statistic
``pom.'' For each of the polynomials $\sum_{k}g_{n}(k)x^k$,
$\sum_{k}h_{n}(k)x^k$, we then have {\it three} combinatorial interpretations: one on alternating permutations by the statistic ``grn,'' and two on labeled, binary, increasing, topological trees by ``pom'' and ``eoc.''

\medskip
\goodbreak

Recently, there has been a revival of studies on arithmetical and combinatorial properties of both tangent and secant numbers.  Ordering the alternating 
permutations according
to their leftmost elements has led to the {\it Entringer recurrence},
having interesting properties ([En66, FZ71, FZ71a, Po97, GHZ10]). The geometry of those permutations has been fully exploited
([KPP94, St10]), in particular by looking at their
quadrant marked mesh patterns in [KR12], or
for defining  and studying natural
$q$-{\it analogs} of the tangent and secant numbers ([AG78, AF80, 
Fo81], [St99, p.~148-149]).  Further
$q$-analogs were also introduced, based no longer on alternating
permutations, but on the so-called {\it doubloon} model (see 
 [FH10, FH10a, FH11]). The classical continued fraction expansions of
secant and tangent have made possible the discovery of other 
$q$-analogs  (see  [Pr08, Pr00, Fu00, HRZ01, Jos10, SZ10]).

\bigskip

{\it Acknowledgement}. The authors should like to thank the referee
for his careful reading and his knowledgeable remarks.

\vglue.6cm
\bigskip\bigskip
\centerline{\bf References}

{\eightpoint

\bigskip

\article AF80|George Andrews; Dominique Foata|Congruences for the
$q$-secant number|Europ. J. Combin.|1|1980|283--287|

\article AG78|George Andrews; Ira Gessel|Divisibility
properties of the $q$-tangent numbers|Proc. Amer. Math.
Soc.|68|1978|380--384|

\article  An1879|D\'esir\'e Andr\'e|D\'eveloppement de $\sec x$ et
$\tan x$|C. R. Math. Acad. Sci. Paris|88|1879|965--979|

\article An1881|D\'esir\'e Andr\'e|Sur les permutations
altern\'ees|J. Math. Pures et Appl.|7|1881|167--184|

\livre Co74|Comtet, Louis|\hskip-5pt Advanced 
Combinatorics|\hskip-5pt D.
Reidel/Dordrecht-Holland, Boston, {\oldstyle 1974}|

\article En66|R.C. Entringer|A combinatorial interpretation of the
Euler and Bernoulli numbers|Nieuw Arch. Wisk.|14|1966|241--246|

\article Fo81|Dominique Foata|\kern -2pt Further
divisibility properties of the $q$-tangent numbers|Proc. Amer. Math.
Soc.|81|1981|143--148|

\article FH10|Dominique Foata; Guo-Niu Han|\kern-2pt The doubloon
polynomial triangle|Ramanujan J.|23|2010|107--126 (The Andrews
Festschrift)|

\article FH10a|Dominique Foata; Guo-Niu Han|Doubloons and
$q$-secant numbers|M\"unster J. of Math.|3|2010|129--150|

\article FH10b|Dominique Foata; Guo-Niu Han|The $q$-tangent and
$q$-secant numbers via basic Eulerian polynomials|Proc.
Amer. Math. Soc.|138|2010|385-393|

\article FH11|Dominique Foata; Guo-Niu Han|Doubloons and new
$q$-tangent numbers|Quarterly J.  Math.|62|2011|417--432|

\divers FS71|Dominique Foata; Marcel-Paul Sch\"utzenberger|Nombres
d'Euler et permutations alternantes. Manuscript (unabridged version)
71 pages, University of Florida, Gainesville, {\oldstyle 1971}.\hfil\break
({\tt
http://igd.univ-lyon1.fr/\char126slc/books/index.html})|

\divers FS71a|Dominique Foata; Marcel-Paul Sch\"utzenberger|Nombres
d'Euler et permutations alternantes, in J. N. Srivastava et {\it al.}
(eds.), {\it A Survey of Combinatorial Theory}, North-Holland,
Amsterdam, 1973, pp. 173-187|

\divers Fu00|Markus  Fulmek|A continued fraction expansion for
a $q$-tangent function, {\sl S\'em. Lothar. Combin.},
{\bf B45b} ({\oldstyle2000}), 3pp|

\article GHZ10|Yoann Gelineau; Heesung Shin; Jiang Zeng|Bijections for
Entringer families|Europ. J. Combin.|32|2011|100--115|

\divers Ha12|Guo-Niu Han|The Poupard Statistics on Tangent and Secant Trees, Strasbourg, preprint 12~p|

\divers HRZ01|Guo-Niu Han; Arthur Randrianarivony; Jiang Zeng|Un
autre $q$-analogue des nombres d'Euler, {\sl The Andrews Festschrift.
Seventeen Papers on Classical Number Theory  and Combinatorics}, D.
Foata, G.-N. Han eds., Springer-Verlag,  Berlin Heidelberg, {\oldstyle
2001}, pp. 139-158. {\sl S\'em. Lothar. Combin.}, Art. B42e, 22 pp|

\livre Jo39|Charles Jordan|Calculus of Finite Differences|R\"ottig and
Romwalter,  Budapest, {\oldstyle 1939}|

\article Jos10|M. Josuat-Verg\`es|A $q$-enumeration of
alternating permutations|Europ. J. Combin.|31|2010|1892--1906|

\divers KR12|Sergey Kitaev; Jeffrey Remmel|Quadrant Marked Mesh Patterns in Alternating Permutations, {\sl S\'em. Lothar. Combin.},
{\bf B68a} ({\oldstyle 2012}), 20pp|

\article KPP94|A. G. Kuznetsov; I. M. Pak;  A. E.
Postnikov|Increasing trees and alternating permutations|Uspekhi Mat.
Nauk|49|1994|79--110|

\livre Ni23|Niels Nielsen|Trait\'e \'el\'ementaire des nombres
de Bernoulli|Paris, Gauthier-Villars, {\oldstyle 1923}|

\article Po82|Christiane Poupard|De nouvelles significations
\'enum\'eratives des nombres\hfil\break d'Entringer|Discrete
Math.|38|1982|265--271|

\article Po89|Christiane Poupard|Deux propri\'et\'es des arbres
binaires ordonn\'es stricts|Europ. J. Combin.|10|1989|369--374|

\article  Po97|Christiane Poupard|Two other interpretations of the
Entringer numbers|Europ. J. Combin.|18|1997|939--943|

\article Pr00|Helmut Prodinger|Combinatorics of geometrically
distributed random variables: new $q$-tangent and $q$-secant
numbers|Int. J. Math. Math. Sci.|24|2000|825--838| 

\divers Pr08|Helmut Prodinger|A Continued Fraction Expansion
for a $q$-Tangent Function: an Elementary Proof, {\sl S\'em.
Lothar. Combin.}, {\bf B60b} ({\oldstyle 2008}), 3 pp|

\article SZ10|Heesung Shin; Jiang Zeng|The $q$-tangent and
$q$-secant numbers via continued fractions|Europ. J. 
Combin.|31|2010|1689--1705|

\divers Sl07|N.J.A. Sloane|On-line Encyclopedia of Integer
Sequences,\hfil\break
{\tt
http://www.research.att.com/\char126 njass/sequences/}|

\livre St99|Richard P.
Stanley|Enumerative Combinatorics, vol. 2|
Cambridge University press, {\oldstyle 1999}|

\divers St10|Richard P.
Stanley|A Survey of Alternating Permutations, in {\sl Combinatorics and graphs}, 165--196, {\sl Contemp. Math.}, {\bf 531}, Amer. Math. Soc. Providence, RI, {\oldstyle 2010}|

\divers Vi88|Xavier G. Viennot|S\'eries g\'en\'eratrices
\'enum\'eratives, chap.~3, Lecture Notes, 160~p., 1988, notes de
cours donn\'es
\`a l'\smash{\'E}cole Normale Sup\'erieure Ulm (Paris), UQAM (Montr\'eal,
Qu\'ebec) et Universit\'e de Wuhan (Chine)\hfil\break
{\tt
http://web.mac.com/xgviennot/Xavier\_Viennot/cours.html}|

\bigskip
\hbox{\vtop{\halign{#\hfil\cr
Dominique Foata \cr
Institut Lothaire\cr
1, rue Murner\cr
F-67000 Strasbourg, France\cr
\noalign{\smallskip}
{\tt foata@unistra.fr}\cr}}
\qquad
\vtop{\halign{#\hfil\cr
Guo-Niu Han\cr
I.R.M.A. UMR 7501\cr
Universit\'e de Strasbourg et CNRS\cr
7, rue Ren\'e-Descartes\cr
F-67084 Strasbourg, France\cr
\noalign{\smallskip}
{\tt guoniu.han@unistra.fr}\cr}}}

}

\bye